\documentclass[12pt]{article}
\usepackage{amsfonts, amssymb, amsmath}

\usepackage{enumerate}

\usepackage[colorlinks=true, urlcolor=blue,linkcolor=blue, citecolor=blue]{hyperref}

\usepackage{graphics,graphicx,amsmath}

\newcommand{\N}{\mathbb{N}}

\newtheorem{prop}{Proposition}[section]

\newtheorem{corollary}[prop]{Corollary}
\newtheorem{theorem}[prop]{Theorem}
\newtheorem{remark}[prop]{Remark}

\def\({\left(}
\def\){\right)}

\def\[{\left[}
\def\]{\right]}

\newcommand{\Z}{\mathbb{Z}}
\newcommand{\B}{\mathbb{B}}

\newcommand{\genstirlingII}[3]{%
  \genfrac{\{}{\}}{0pt}{#1}{#2}{#3}%
}
\renewcommand{\S}[2]{\genstirlingII{}{#1}{#2}}

\newcommand{\genstirlingI}[3]{%
  \genfrac{[}{]}{0pt}{#1}{#2}{#3}%
}
\newcommand{\s}[2]{\genstirlingI{}{#1}{#2}}

\renewcommand{\=}{\equiv}

\newenvironment{Proof}{\removelastskip\par\medskip
\noindent{\em Proof.} \rm}{\penalty-20\null\hfill$\square$\par\medbreak}

\allowdisplaybreaks

\numberwithin{equation}{section}

\begin{document}

\title{
\huge
Backward Touchard congruence
} 

\author{
Grzegorz Serafin\thanks{Faculty of Pure and Applied Mathematics, Wroc{\l}aw University of Science and Technology, Ul. Wybrze\.ze Wyspia\'nskiego 27, Wroc{\l}aw, Poland.
e-mail: {\tt grzegorz.serafin@pwr.edu.pl}}}

\maketitle

\begin{abstract} 
The celebrated Touchard congruence states that $B_{n+p}\=B_n+B_{n+1}$ modulo $p$, where $p$ is a prime number and  $B_n$ denotes the Bell number. In this paper we study divisibility properties of $B_{n-p}$ and their generalizations involving  higher powers of $p$ as well as the $r$-Bell numbers. In particular, we show a closely relation of the considered problem  to the Sun-Zagier congruence, which is additionally improved by deriving \mbox{a new} relation between $r$-Bell and derangement numbers. Finally, we conclude some results on the  period of the Bell numbers modulo $p$. 
\end{abstract}

\baselineskip0.7cm

\section{Introduction}
The Bell numbers $B_n$ are one of the most classical sequences in Combinatorics and describe the number of partitions of a given set of $n$ elements into non-empty subsets.  Their various aspects have been studied for more than a hundred years. In particular, the divisibility properties are of special interest (see, among others,  \cite{GR, G, Kahale, LPS, MR, SWZ, SZ, Touchard, Tsumura}). The first remarkable result in this direction comes from 1933, when Jackues Touchard \cite{Touchard} obtained the congruence 
\begin{align}\label{eq:Touchard}
B_{n+p^m}\=mB_n+B_{n+1}\ \ \ \ \ (\mbox{mod }p), 
\end{align}
valid for any natural $n,m$ and prime $p$.  Nevertheless, the term 'Touchard congruence' refers usually to the case $m=1$ only. Another interesting relation was discovered by Sun and Zagier  \cite{SZ}
\begin{align}\label{eq:congr1}
\sum_{k=1}^{p-1} \frac{B_k}{(-m)^k}\= (-1)^{m-1}D_{m-1}, \ \ \ \ \ \ \ \ \ (\mbox{mod }p\Z_p),
\end{align}
where $m\geq1$, $p\nmid m$ is a prime number, $D_n$ stands for the $n$-th derangement number and  $\Z_p$ denotes the ring of $p$-adic integers.

The $r$-Bell numbers $B_{n,r}$, $r\geq0$,  are a natural generalization of the Bell numbers and count partitions of a set of $n+r$ elements such that $r$ chosen elements are separated \cite{M1}. 
The case $r=0$ clearly corresponds to the standard Bell numbers $B_{n,0}=B_{n}$.  It turns out thatthe $r$-Bell numbers satisfy \eqref{eq:Touchard} as well \cite{MR}.  
Furthermore, the Sun-Zagier congruence \eqref{eq:congr1} has been  improved \cite{BenM,MR,Sun,SWZ} reaching the following form involving the $r$-Bell numbers \cite{Serafin}
\begin{align}\label{eq:SZS}
\sum_{i=1}^{p^a-1}\frac{B_{n+i,r}}{(-m)^i}\=a\sum_{k=0}^n\S {n}{k}_r(-1)^{k+m+r-1}D_{k+m+r-1},\ \ \ \ \ (\mbox{mod}\  p\Z_p ),
\end{align}
where $a\geq1$ and $\S {n}{k}_r$ are the $r$-Stirling numbers of the second kind - see Section \ref{sec:Stirling} for more details. In fact, these results have been obtained in the polynomial version as well, where $r$-Bell and derangement numbers are replaced by $r$-Bell and derangement polynomials, respectively. Note that the $r$-Bell numbers appear naturally when studying divisibility properties of the classical Bell numbers, which has been shown in \cite{Serafin} and will be confirmed in \mbox{Section \ref{sec:btc}}.

In this article, we study divisibility properties of $B_{n-p^m}$, which might be described  as backward analog of  \eqref{eq:Touchard}, as well as their generalizations for the $r$-Bell numbers. In particular, we show their relation to the generalized Sun-Zagier congruence \eqref{eq:SZS}, which is additionally improved by providing equivalent forms where number of terms in the sum does not depend on $n$. This is achieved by deriving some new identities bonding $r$-Bell and derangement numbers. 

One of the simplest and, at the same time, most elegant results of the paper (see Corollary \ref{cor:btc}) states that
$$B_{n-p}\= V_n,\ \ \ \ (\mbox{mod }p),$$
where $V_n$ is the number of partitions of the set $\{1,2, ...,n\}$ without singletons. Apparently, $V_n$ may be considered as $B_{n,-1}$. It turns out to be a part of a more general rule, which motivated us to introduce the $r$-Bell numbers for negative values of the index $r$. In fact, this could be avoided due to the relation $B_{n,r}\=B_{n,r+p}$ (mod $p$), $r\in\Z$, however, it is definitely clearer that $B_{n,-l}$ is independent of $p$ than $B_{n,p-l}$ is, for some $l\geq0$. Furthermore, we apply the established equivalences to conclude some result on periodicity of the sequences $B_{n,r}$ \mbox{mod $p$.} In particular, we address the hypothesis that  $\frac{p^p-1}{p-1}$ is their minimal period by excluding a class of some other potential periods.

\section{Preliminaries}
\subsection{The $r$-Stirling numbers}
\label{sec:Stirling}
We denote by $\s nk_r$ and $\S nk_r$, $r\geq0$, the $r$-Stirling numbers of the first and second kind, respectively. The number $\s nk_r$ counts permutations of the set $\{1,2,...,n+r\}$ having $k+r$ cycles such that the numbers $1,2,...,r$ are in distinct cycles, while  $\S nk_r$ represents  the number of partitions of the set $\{ 1, 2,  . . .  , n+r\}$ into $k+r$
non-empty disjoint subsets, such that the numbers $1, 2 , . . . , r$ are in distinct subsets. For $r=0$ we obtain the classical Stirling numbers. The $r$-Stirling numbers were introduced and described in details by A. Z. Broder in \cite{B} (see also \cite{CI, CII}), where slightly different notation was used ($\s {n+r}{k+r}_r$ instead of $\s nk_r$ and similarly $\S {n+r}{k+r}_r$ instead of $\S nk_r$). Nevertheless, the convention used in this paper seems more natural and makes expressions less complicated.

The $r$-Stirling numbers may be also characterized by the following expansions 
\begin{align}\label{eq:(x+r)^n}
(x+r)^n&=\sum_{k=0}^n\S nk_r x^{\underline k},\\
\label{eq:(x+r)^_n}
(x+r)^{\overline n}&=\sum_{k=0}^n \s{n}{k}  _r x^{ k},
\end{align}
where $x^{\underline k}=x(x-1)\cdot...\cdot(x-k+1)$ and $x^{\underline k}=x(x+1)\cdot...\cdot(x+k-1)$ are the falling and rising factorials, respectively. They admit the orthogonality relation  
\begin{align*}
\sum_{k=m}^n\s nk_r\S km(-1)^k=\sum_{k=m}^n\S nk_r\s km(-1)^k=(-1)^n\delta_{mn}.
\end{align*}
The exponential generating functions are given by
\begin{align*}
\sum_{n=k}^\infty \s nk_r\frac{x^n}{n!}&=\frac1{k!}\(\frac1{1-x}\)^r\(\ln\(\frac1{1-z}\)\)^k,\\
\sum_{n=k}^\infty \S nk_r\frac{x^n}{n!}&=\frac1{k!}e^{rx}\(e^x-1\)^k.
\end{align*}
Consequently, defining the $r$-Stirling transform (involving the $r$-Stirling numbers of the second kind only)
 of a sequence $(a_n)_{n\geq1}$ by
 $$b_n=\sum_{k=1}^n\S{n}{k}_ra_k,$$
there holds the following relation between the  exponential generating functions $A(x)$, $B(x)$ of $(a_n)$ and $(b_n)$, respectively,
\begin{align*}
B(x):=\sum_{n=1}^\infty \frac{b_n}{n!}x^n=e^{rx}A(e^x-1).
\end{align*}


\subsection{The $r$-Bell numbers }

As mentioned in Introduction, the $r$-Bell numbers $B_{n,r}$, $r\geq0$, count partitions of the set $\{1,...,n+r\}$ into  
non-empty disjoint subsets such that the numbers $1,...,r$ are separated. In particular, we  have
$$B_{n,r}=\sum_{i=0}^n\S{n}{k}_r.$$
The exponential generating function takes the form
\begin{align}\label{eq:egfB_r}
\mathcal B_r(t)=\sum_{n=0}^\infty  B_{n,r}\frac{t^n}{n!}=e^{e^t-1+rt}.
\end{align}
Treating it as the definition, we can clearly  extend the range of $r$ onto all integer numbers $\Z$. The combinatorial interpretation of the Bell numbers with a negative index $r$ is only known, up to the authors knowledge, in the case $r=-1$. Namely, $B_{n,-1}$, denoted usually by $V_n$, represents the number of partitions of the set of $n$ elements containing no singletons.

Applying the general Leibniz rule to \eqref{eq:egfB_r}, we obtain
\begin{align}\label{eq:B=sumbB}
\sum_{k=0}^{n}{n \choose k}B_{k,r} \,m^{n-k}=B_{n,r+m}
\end{align}
for all $n\geq0$, $r\in\Z$. As a consequence, we get the following periodicity property  
\begin{align}\label{eq:r+p}
B_{n,r+p}\= B_{n,r}, \ \ \ \ (\mbox{mod }p),
\end{align}
valid for any prime $p$ and $n\geq0$, $r\in \Z$. In particular, applying this to Theorem 5 in \cite{MR}, which generalized \eqref{eq:Touchard} onto $r$-Bell numbers for $r\geq0$, we arrive at
\begin{align}\label{eq:r-Touchard}
B_{n+p^a,r}\=B_{n+1,r}+aB_{n,r}, \ \ \ \ \ \ \ \ \ (\mbox{mod}\  p),
\end{align}
where $a,n\in\N$ and $r\in\Z$ is any intiger number. We close this section with another useful recurrences (\cite{CI}, eq. (3.22-3.23)) 
\begin{align}\label{eq:Bn+m,r}
B_{n+m,r}&=\sum_{k=0}^m\S {m}{k}_r x^{k}B_{n,k+r},\\\label{eq:Bn,m+r}
B_{n,r+m}&=\sum_{k=0}^m(-1)^{m-k}\s {m}{k}_r x^{k}B_{n+k,r},
\end{align}
where $m,n,r\geq0$.

\section{Generalized Sun-Zagier congruence}
 The main feature of  the generalized Sun-Zagier congruence \eqref{eq:SZS} is that the right-hand side  does not depend on $p$. Nevertheless, since $r$ and $m$ are typically supposed to be fixed and $n$ may vary, the number of terms in the sum could be arbitrarily large.  In the next theorem we   present two identites removing this inaccuracy.

\begin{theorem}\label{thm:sumD}
For $n \geq0$ and $r+m\geq1$ we have 
\begin{align*}
\sum_{k=0}^n\S {n}{k}_r(-1)^{k+m+r-1}D_{k+m+r-1}&=\sum_{k=0}^m\sum_{i=0}^r\s mk\s ri(-1)^{r-i}B_{n+k+i-1,-m}\\
&=\sum_{k=0}^{m+r-1}(-1)^kk!B_{n,-k-1},\ \ \ \ \ \ \ \ \ \  (\mbox{mod}\  p).
\end{align*}
\end{theorem}

\begin{Proof}
In order to derive the first identity from the assertion we will employ the umbral calculus. It allows us to represent $B_n$ as $\B^n$, where $\B$ is a symbol called umbra. Such an approach has been already used in e.g. in \cite{BenM, MR,SW,SWZ} in the context of Bell numbers. 

Taking $r=0$ in \eqref{eq:B=sumbB} and by the binomial theorem we obtain
\begin{align}\label{aux1}B_{n,m}=\(\B+m\)^n,\ \ \ \ \ \ \ \ \ m\in\Z.\end{align}
Furthermore, from Lemma 2.2 in \cite{SWZ} 
we know that for any $k\geq0$ it holds
\begin{align}\label{eq:B-1}
(\B-1)^{\underline{k}}=(\B-1)(\B-2)...(\B-k)=(-1)^kD_k.
\end{align}
 Thus,  we get 
\begin{align*}
&\sum_{k=0}^n\S {n}{k}_r(-1)^{k+r+m-1}D_{k+r+m-1}\\
&=\sum_{k=0}^n\S {n}{k}_r(\B-1)(\B-2)...(\B-k-r-m+1)\\
&=(\B-1)^{\underline{r+m-1}}\sum_{k=0}^n\S {n}{k}_r(\B-r-m)^{\underline k}\\
&=(\B-1)^{\underline{r+m-1}}(\B-m)^n,
\end{align*}
where we used \eqref{eq:(x+r)^n}. For $n\geq1$, we rewrite it  as follows
\begin{align*}
(\B-1)^{\underline{r+m-1}}(\B-m)^n&=\big((\B-m)+(m-1)\big)^{\underline{r+m-1}}(\B-m)^n\\
&=\big((\B-m)+(m-1)\big)...\big((\B-m)-r+1\big)(\B-m)^n\\
&=(\B-m)^{\overline m}(-1)^r(m-\B)^{\overline r}(\B-m)^{n-1}.
\end{align*}
 Hence,  the identity   \eqref{eq:(x+r)^_n}   gives us
\begin{align*}
L&=\sum_{k=0}^m\s mk(\B-m)^{k}(-1)^r\sum_{i=0}^r\s ri(m-\B)^{i}(\B-m)^{n+k+i-1}\\
&=\sum_{k=0}^m\sum_{i=0}^r\s mk\s ri(-1)^{r-i}(\B-m)^{n+k+i-1}\\
&=\sum_{k=0}^m\sum_{i=0}^r\s mk\s ri(-1)^{r-i}B_{n+k+i-1,-m},
\end{align*}
as required. Eventually, let us observe that due to the assumption $r+m\geq0$ the above calculations hold true for $n=0$. Indeed, assuming $m\geq1$, we  have $\s m0=0$ and  $x^{\overline m}=xP(x)$ for some polynomial $P$, so from \eqref{eq:(x+r)^_n} we have
$$x^{\overline m}x^{-1}=\sum_{k=0}^m \s{m}{k}  _r x^{ k-1},$$
which may be applied to $(\B-m)^{\overline m}(\B-m)^{-1}$. We proceed similarly in the case $r\geq1$, $m=0$.

To obtain the other identity from the assertion, let us write
\begin{align}\label{aux5}\sum_{k=0}^n\S {n}{k}_r(-1)^{k+m+r-1}D_{k+m+r-1}=\sum_{k=0}^n\S {n}{k}_ra_k,\end{align}
where 
$$a_k=(-1)^{k+m+r-1}D_{k+m+r-1}.$$
Since $e^{-t}/(1-t)$ is the exponential generating function of the sequence $D_k$,  the exponential generating function of $a_k$ is given by
\begin{align*}
A(x)&=\frac{d^{m+r-1}}{dx^{m+r-1}}\frac{e^t}{1+t}=\sum_{i=0}^{m+r-1}\frac{e^t}{(1+t)^{1+i}}i!(-1)^i,
\end{align*}
where we used the general Leibniz rule. Hence, since \eqref{aux5} is the $r$-Stirling transform of $a_k$, its exponential generating function takes the  form
$$e^{rt}A(e^t-1)=\sum_{k=0}^{r+m-1}\frac{e^{e^t-1}}{(e^t)^{1+k}}k!(-1)^k=\sum_{k=0}^{r+m-1}(-1)^kk!{e^{e^t-1-(k+1)t}}.$$
We can identify the exponents in the sum as exponential generating functions of the $r$-Bell numbers \eqref{eq:egfB_r}, which ends the proof.
\end{Proof}
\begin{remark}
Repeating the arguments from the proof, one can easily obtain a 'polynomial' version of the theorem, i.e. with $B_{n,r}(x)=\sum_{i=0}^n\S{n}{k}_rx^k$ and $D_n(x)=\sum_{k=0}^n{n \choose k}k!(x-1)^{n-k}$ instead of $B_{n,r}$ and $D_{n}$, respectively. It is enough to know that the corresponding exponential generating functions are $e^{x(e^t-1)+rx}$ and $e^{(x-1)t}/(1-t)$.
\end{remark}

The first identity from Theorem \ref{thm:sumD} takes especially simple form for $m=0$. Namely, by virtue of \eqref{eq:Bn,m+r}, we get
\begin{corollary}\label{cor:1}
For $n \geq1$ and $r\geq1$ it holds
\begin{align*}
\sum_{k=0}^n\S {n}{k}_r(-1)^{k+r-1}D_{k+r-1}=B_{n-1,r}.
\end{align*}
\end{corollary}

\begin{remark}\label{rem:1} The above  result is valid for $r=0$ as well. This very special case, in the 'polynomial' version,  is covered by Lemma 2.2 in \cite{SWZ}. 
\end{remark}

Finally, let us formulate the new version of the Sun-Zagier congruence.  Applying Theorem \ref{thm:sumD} to \eqref{eq:SZS}, we obtain
\begin{corollary} For any  $a,m\geq 1$, $n\geq0$ and any prime number $p \nmid m$,
we have
\begin{align*}
\sum_{i=1}^{p^a-1}\frac{B_{n+i,r}}{(-m)^i}&\=a\sum_{k=0}^m\sum_{i=0}^r\s mk\s ri(-1)^{r-i}B_{n+k+i-1,-m}\\
&\=a\sum_{k=0}^{r+m-1}(-1)^kk!B_{n,-k-1},\ \ \ \ \ \ \ \ \ \ \ \ \ \ \ \ \ \ (\mbox{mod}\  p\Z_p).
\end{align*}
\end{corollary}

\section{Backward Touchard congruence}
\label{sec:btc}
We start this section with an equivalence, which simplifies \eqref{eq:r-Touchard} in some cases and was one of the motivations of research presented in the article.
\begin{prop}
For $n,r\geq0$ we have 
$$B_{n+p^r,-r}\=B_{n,-r+1}, \ \ \ \ \ \ \ \ \ (\mbox{mod}\  p).$$
\end{prop}
\begin{Proof}
From \eqref{eq:Bn+m,r} with $m=1$ we get for any $l\in\Z$
$$B_{n,l+1}\=B_{n+1,l}-lB_{n,l}, \ \ \ \ \ \ \ \ \ (\mbox{mod}\  p).$$
Hence, substituting $l=-r\leq0$ and using \eqref{eq:r-Touchard}, we get
$$B_{n,-r+1}\=B_{n+1,-r}+rB_{n,-r}\=B_{n+p^r,-r}, \ \ \ \ \ \ \ \ \ (\mbox{mod}\  p),$$
which ends the proof.
\end{Proof}

In particular, for $r=1$ we obtain the below-given elegant congruence, which may be called the backward Touchard congruence.
\begin{corollary}\label{cor:btc}
For a prime $p$ and natural $n\geq p$ we have
$$B_{n-p}\=V_n, \ \ \ \ \ \ \ \ \ (\mbox{mod}\  p).$$

\end{corollary}
 Next, using simple induction argument, one can generalize it as follows.

\begin{corollary} \label{prop:sum}
We have
$$B_{n-\sum_{k=1}^rp^k}\=B_{n,-r}, \ \ \ \ \ \ \ \ \ (\mbox{mod}\  p).$$
\end{corollary}
This equivalence has further consequences on the period of the $r$-Bell numbers modulo $p$ - see the next section for details.
Nevertheless,  the most natural direction of research is to investigate divisibility properties of  $B_{n-p^m,r}$, which is executed in the next theorem.
\begin{theorem}\label{thm:bT=D}
For a prime $p$ and  integers $n,m\geq1$, $r\geq0$ such that $n\geq p^m$ and $p\nmid m$ we have
\begin{align*}
B_{n-p^m,r}&\=\sum_{k=0}^n\S {n}{k}_r(-1)^{k+m+r-1}D_{k+m+r-1}\\
&\=\sum_{k=0}^m\sum_{i=0}^r\s mk\s ri(-1)^{r-i}B_{n+k+i-1,-m}\\
&\= \sum_{k=0}^{r+m-1}(-1)^kk!B_{n,-k-1},\ \ \ \ \ \ \ \ \ \ \ \ \ \ \ \ \ \ \ \ \ \  (\mbox{mod}\  p).
\end{align*}
\end{theorem}
\begin{Proof}
We will show only the first equivalence. The other ones follow from Theorem \ref{thm:sumD}. 

From \eqref{eq:r-Touchard} we have
\begin{align}\label{aux2}mB_{n,r}\=B_{n+p^m,r}-B_{n+1,r}, \ \ \ \ \ \ \ \ \ (\mbox{mod}\  p).\end{align}
More generaly, for any $N\geq 1$ it holds
\begin{align}\label{aux3}
m^NB_{n,r}\=\(\sum_{k=0}^{N-1}(-1)^km^{N-1-k}B_{n+p^m+k,r}\)+(-1)^NB_{n+N,r}, \ \ \ \ \ \ \ \ \ (\mbox{mod}\  p),
\end{align}
which may be shown by induction. Indeed, multiplying it by $m$ and using \eqref{aux2}, we get
\begin{align*}
m^{N+1}B_{n,r}&\=\(\sum_{k=0}^{N-1}(-1)^km^{N-k}B_{n+p^m+k,r}\)+(-1)^NmB_{n+N,r}\ \ \ \ \ \ \ \ \ \ \ \ \ \ \ \ \  (\mbox{mod}\  p)\\
&=\(\sum_{k=0}^{N-1}(-1)^km^{(N+1)-1-k}B_{n+p^m+k,r}\)+(-1)^N\(B_{n+p^m+N,r}-B_{n+N+1,r}\)\\
&=\(\sum_{k=0}^{(N+1)-1}(-1)^km^{(N+1)-1-k}B_{n+p^m+k,r}\)+(-1)^{N+1}B_{n+(N+1),r},
\end{align*}
as required. Next, substituting $N=p^m$ in \eqref{aux3} and exploiting the congruence $m^{p^m}\= m$, valid by virtue of   Fermat's little theorem, we arrive at
\begin{align*}
mB_{n,r}\=m^{p^{m}}B_{n,r}&\=\(\sum_{k=0}^{p^m-1}(-1)^km^{p^m-1-k}B_{n+p^m+k,r}\)+(-1)^{p^m}B_{n+p^m}\\
&\=\sum_{k=1}^{p^m-1}\frac{B_{n+p^m+k,r}}{(-m)^k}+\big(1+(-1)^{p^m}\big)B_{n+p^m+k,r}, \ \ \ \ \ \ \ \ \ \ \ \ \ \ \ (\mbox{mod}\  p\Z_p).
\end{align*}
The equivalence $1+(-1)^{p^m}\=0$ mod $p$, which holds for any prime $p$, and  application of the generalized Sun-Zagier congruence \eqref{eq:SZS} end the proof.
\end{Proof}

\begin{remark} 
For $m=0$ we the first equivalence from Theorem \ref{thm:bT=D} becomes an identity  - see Corollary \ref{cor:1} and Remark \ref{rem:1}. However,  it is not true anymore for $m\geq1$. For example, one can verify that for $m=1$, $r=0$, $p=2$, $n=4$.
\end{remark}


\section{The $r$-Bell numbers modulo a prime number}
In this section we focus on the periodicity of the sequence $B_{n,r}$ mod $p$, for a prime $p$, which is simply the sequence of reminders of a division $B_{n,r}$ by $p$. Hall \cite{Hall} discovered that $B_n$ has the period
$$N_p=\frac{p^p-1}{p-1}=1+p+p^2+...+p^{p-1}.$$
We can easily recover it combining Corollary \ref{prop:sum} with \eqref{eq:r+p} and the equality $B_{n,1}=B_{n+1}$:
$$B_{n-N_p}=B_{n-1-\sum_{k=1}^{p-1}p^k}\=B_{n-1,-p+1}\=B_{n-1,1}=B_n, \ \ \ \ \ \ \ \ \ (\mbox{mod}\  p).$$
  Mez\H{o} and Ram\'{\i}rez \cite{MR} extended the Hall's result for $r$-Bell numbers by showing that $N_p$ is a period of $B_{n,r}$, for a fixed $r\geq 0$.  Apparently,  the relation between 
the sequences for different $r$'s is much stronger, which is presented below.
 
\begin{corollary}
For any $r\in\Z$ the sequence $B_{n,r}$ mod $p$ is equal to the sequence $B_n$ mod $p$ shifted in the following manner
$$B_{n,r}\=B_{n-K}, \ \ \ \ \ \ \ \ \ (\mbox{mod}\  p),\ \ \ \ \ n\geq K,$$
where $K=\sum_{k=1}^{(-r \text{ mod } p)}p^k$.
\end{corollary}
\begin{Proof}
 The case $r\in\{0,...,p-1\}$ is cover by Corollary \ref{prop:sum}. One can extend it onto all $r\in\Z$ by virtue of \eqref{eq:r+p}.
\end{Proof}

$N_p$ was proven \cite{LD, MNW,W} to be the minimal period for $p<126$ as well as for $p=137,149,157,163,167$ and $173$. For other primes the problem is open, however, there exist in the literature some partial results. They are usually related to the divisibility properties of $N_p$. A quantitative bound  \cite{LPS,CGRV} states that the minimal period is greater than 
$$\frac12{2p \choose p}+p.$$
In the next theorem, we provide a result of new type, referring to the representation of the period in the base $p$ numerical system.

\begin{theorem}
Let 
\begin{align}\label{aux4}
P_p=\sum_{k=0}^{p-2}a_kp^k<p^{p-1},\ \ \  a_k\in\{0,1,...,p-1\},
\end{align}
 be a  period (not necessarily the  minimal one) of  $B_n$ modulo a prime $p$. Then 
$$\sum_{k=0}^{p-2}a_k\geq p+1.$$ 
\end{theorem}
\begin{Proof}First, we will justify that we can assume $a_0>0$. Namely, if $a_0=a_1=...=a_{k_0}=0$ and $a_{k_0+1}>0$, then $P_p/p^{k_0}$ is a period as well with the same sum of digits (in the base $p$ numerical system). This follows from the fact that both: $N_p$ and $P_p$ are multiplicities
 of the minimal period, while  all the dividers of $N_p$ are of the form $2kp+1$, $k\in\N$ (see \cite{D}, p. 381), and hence the minimal period is not divisible by $p$, so $P_p/p^{k_0}$ has to a multiplicity of the minimal period.
 
  Denote $M=\sum_{k=1}^{p-2}a_k.$
Next, starting from $B_{n+P_p}$, we exploit  the congruence \eqref{eq:Touchard}
$$B_{n+p^m}\=B_{n+1}+mB_n, \ \ \ \ \ \ \ \ \ (\mbox{mod}\  p),$$
in $M-a_0$ steps. In each step we apply it to every term that already appeared and with the same $m\geq1$. Eventually, this procedure leads to
\begin{align}\label{eq:B-B}
B_{n+P_p}\=\sum_{k=a_0}^{M}b_kB_{n+k}, \ \ \ \ \ \ \ \ \ (\mbox{mod}\  p),
\end{align}
for some $b_k\geq0$, $k\in\{0,...,M\}$, such that $b_{M}=1$. Since the sum in \eqref{aux4} is up to $p-2$,  in each step the sum of coefficients of appearing Bell numbers increases $(m+1)$ times for some $2\leq m\leq p-1$. Thus,  the sum $S:=\sum_{k=a_0}^{M}b_k$ is a product of $M-a_0$ positive numbers smaller than $p$, and consequently $S\not\equiv 0$ mod $p$.
Furthermore, since $P_p$ is a period, we get
$$\sum_{k=a_0}^{M}b_kB_{n+k}\=B_n, \ \ \ \ \ \ \ \ \ (\mbox{mod}\  p).$$
Let as suppose that $M\leq p$. Subtracting $B_n$ from both sides and applying the Touchard congruence to $B_{n+M}$, we obtain for $n\geq p-M$
$$\sum_{k=M-p}^{M-1}c_kB_{n+k}\=0, \ \ \ \ \ \ \ \ \ (\mbox{mod}\  p),$$
where 
$$\sum_{k=M-p}^{M-1}c_{k}=\sum_{k=a_0}^{M-1}b_k+2b_{M}-1=\sum_{k=a_0}^{M-1}b_k+1=S\not\equiv 0, \ \ \ \ \ \ \ \ \ (\mbox{mod}\  p).$$
This contradicts  Theorem (3.3) in \cite{LPS}, which implies that if $\sum_{k=0}^A\alpha_kB_{n+k}\=0$ for all $n\geq0$, where $A<p$, then $a_k\=0$ for all $k\in\{0,...,A\}$. Putting $\alpha_k=c_{k+M-p}$, $k\in\{0,...,p-1\}$, we have $A=p-1<p$, while $\sum_{k=0}^A\alpha_k\not\equiv 0$ and hence there exists $k\in\{0,...,p-1\}$ such that $\alpha_k\not\equiv0$.
\end{Proof}

\section*{Acknowledgements}
The author was supported by the National Science Centre grant no. 2015/18/E/ST1/00239.
\footnotesize

\def\cprime{$'$} \def\polhk#1{\setbox0=\hbox{#1}{\ooalign{\hidewidth
  \lower1.5ex\hbox{`}\hidewidth\crcr\unhbox0}}}
  \def\polhk#1{\setbox0=\hbox{#1}{\ooalign{\hidewidth
  \lower1.5ex\hbox{`}\hidewidth\crcr\unhbox0}}} \def\cprime{$'$}

\end{document}